\documentclass{amsart}
\usepackage{latexsym,amsxtra,amscd,ifthen}
\usepackage{amsfonts}
\usepackage{verbatim}
\usepackage{amsmath}
\usepackage{amsthm}
\usepackage{amssymb}

\numberwithin{equation}{section}
\theoremstyle{plain}
\newtheorem{theorem}[equation]{Theorem}

\newtheorem{lemma}[equation]{Lemma}

\newtheorem{corollary}[equation]{Corollary}
\newtheorem{proposition}[equation]{Proposition}

\newtheorem{question}[equation]{Question}

\theoremstyle{definition}
\newtheorem{definition}[equation]{Definition}

\newtheorem{remark}[equation]{Remark}

\newtheorem{hypothesis}[equation]{Hypothesis}
\newtheorem{example}[equation]{Example}

\newcommand{\reg}{\operatorname{reg}}
\DeclareMathOperator{\GK}{GKdim}

\DeclareMathOperator{\HB}{H}

\DeclareMathOperator{\rann}{r.ann}

\DeclareMathOperator{\proj}{proj}

\DeclareMathOperator{\Num}{Num}

\DeclareMathOperator{\Hom}{Hom}

\DeclareMathOperator{\coh}{coh}

\DeclareMathOperator{\aut}{Aut}

\newcommand{\wt}{\widetilde}

\newcommand{\cal}{\mathcal}
\newcommand{\mc}{\mathcal}
\newcommand{\mb}{\mathbb}

\DeclareMathOperator{\spec}{Spec}
\newcommand{\lra}{\longrightarrow}

\DeclareMathOperator{\GKtr}{GKtr}

\begin{document}

\begin{abstract}
If $A$ is a strongly noetherian graded algebra generated
in degree one, then there is a canonically constructed graded ring
homomorphism from $A$ to a twisted homogeneous coordinate ring
$B(X, \mc{L}, \sigma)$, which is surjective in large degree.
This result is a key step in the study of projectively
simple rings.
The proof relies on some results concerning the growth of
graded rings which are of independent interest.
\end{abstract}

\title{Canonical maps to twisted rings}

\author{D. Rogalski and J. J. Zhang}

\address{(Rogalski) Department of Mathematics, MIT,
Cambridge, MA 02139-4307, USA}
\thanks{D. Rogalski was partially supported by NSF grant DMS-0202479}

\email{rogalski@math.mit.edu}

\address{(Zhang) Department of Mathematics, Box 354350, University
of Washington, Seattle, WA 98195, USA}
\thanks{J. J. Zhang was partially supported by NSF grant DMS-0245420
and Leverhulme Research Interchange Grant F/00158/X (UK)}

\email{zhang@math.washington.edu}

\subjclass[2000]{16P90, 16S38, 16W50, 14A22}

\keywords{Graded ring, noncommutative projective geometry, twisted
homogeneous coordinate ring, Gelfand-Kirillov dimension, strongly
noetherian}

\maketitle

\section{Introduction}
\label{sect1}

\emph{Twisted homogeneous coordinate rings} (or \emph{twisted rings}
for short) are special noncommutative graded rings, formed from geometric data, 
which were introduced by Artin and Van den Bergh in \cite{AV}.
Given a projective scheme $X$ over a base field $k$, an automorphism
$\sigma$ of $X$ and an invertible sheaf $\mc{L}$ on $X$, 
one defines the twisted ring  
\[
B(X, \mc{L}, \sigma) = \bigoplus_{n = 0}^{\infty} 
\HB^0(X, \mc{L} \otimes \sigma^* \mc{L} 
\otimes \dots \otimes (\sigma^{n-1})^* \mc{L})
\]
which has a natural $\mb{N}$-graded $k$-algebra structure.
If the sheaf $\mc{L}$ is \emph{$\sigma$-ample} (see Section~\ref{sect3} below), 
then the ring $B(X, \mc{L}, \sigma)$
is noetherian and has many other good properties \cite{AV,Ke1}.

Twisted rings actually first appeared 
in \cite{ATV} (though not named as such) as an important
ingredient in the project to classify Artin-Schelter regular
algebras of global dimension 3,
which leads to a classification of quantum projective planes.
Given such a graded regular ring $A$, Artin, Tate and Van den Bergh
showed in \cite{ATV} that there exists a surjective ring homomorphism
$\varphi: A \to B(X, \mc{L}, \sigma)$ for a naturally constructed
triple $(X, \mc{L}, \sigma)$. This is very useful since the
ring $B(X, \mc{L}, \sigma)$, being geometrically defined, may be
studied with geometric techniques.  Then information may pulled
back to the ring $A$. A similar method has been applied with
success to certain classes of graded regular rings of global dimension
4 by several authors \cite{Va, SV, VV1, VV2}.

The goal of this article is to show that 
the method described above also applies to many other 
graded rings.  We will prove that under quite general conditions
on a graded ring $A$, there exists a canonically defined graded ring
homomorphism $\varphi: A \to B(X,\mc{L}, \sigma)$, and that this map
is always surjective in large degree. 
As in \cite{ATV}, the scheme $X$ is constructed as the
solution to a certain moduli problem.  Let $A$ be a graded algebra
generated in degree $1$.  For each commutative algebra $R$, an
$R$-\emph{point module} for $A$ is a cyclic graded right
$A \otimes R$-module
$M = \bigoplus_{n = 0}^{\infty} M_n$ with $M_0 = R$ and such that
$M_n$ is a locally free $R$-module of rank $1$ for all $n \geq 1$.
One may define a functor $F$ from the category of commutative
algebras to the category of sets, by sending each commutative
algebra $R$ to the set of isomorphism classes of $R$-point modules
for $A$.  If the functor $F$ is represented by a projective scheme
$X$, we call $X$ the \emph{point scheme} for the ring $A$. By work
of Artin and Zhang \cite{AZ2} it is known that the point scheme
for a graded ring $A$ exists as long as $A$ is \emph{strongly
noetherian}, which means that $A \otimes_k R$ is noetherian for
all noetherian commutative $k$-algebras $R$.  This last condition
is satisfied for many important classes of noncommutative
graded rings \cite{ASZ}.  Now we may state our main result. 
Except where otherwise noted, in all of the theorems in this paper 
the base field $k$ is
assumed to be algebraically closed.

\begin{theorem}
\label{xxthm1.1}
Suppose $A$ is a strongly noetherian 
connected graded $k$-algebra generated in degree $1$.
\begin{enumerate}
\item
There exists a canonical graded ring homomorphism
$$\varphi: A \to B(X, \mc{L}, \sigma),$$ 
where $X$ is the point scheme
for $A$.  The kernel of $\varphi$ is equal in large degree to the ideal
\[
\{x \in A | Mx = 0\ \text{for all}\ R\text{-point modules}\ M,
\text{all commutative algebras}\ R \}.
\]
\item
Further, the sheaf $\mc{L}$ is
$\sigma$-ample, and $\varphi$ is surjective in large degree.
\end{enumerate}
\end{theorem}

A special case of the above theorem occurs when $A$ is a
graded domain of GK-dimension 2 which is generated in degree 1.  In this case 
the canonical map $\varphi$ has already been constructed in \cite{AS} 
(although in a different way), and $\varphi$ is known to be an isomorphism in large degree 
\cite[Theorem 0.2]{AS}.  
In fact, such rings $A$ are automatically strongly noetherian \cite[Theorem 4.24]{ASZ}, 
and have graded quotient rings 
of the form $Q = K[t^{\pm 1}; \tau]$ for some field $K$ of transcendence
degree $1$ \cite[Theorem 0.1]{AS}.  
Thus the following consequence of our main theorem may be viewed as a generalization of 
\cite[Theorem 0.2]{AS} to higher dimensions.
\begin{theorem}
\label{main cor}
Let $A$ be semiprime, connected graded, generated in degree $1$, and strongly noetherian.  
Suppose that the graded ring of fractions of $A$ has the form $R[t^{\pm 1}; \tau]$, where $R$ is 
a commutative ring.  Then the canonical map $\varphi$ of Theorem~\ref{xxthm1.1} is 
an isomorphism in large degree. 
\end{theorem}

We call a locally finite $\mb{N}$-graded $k$-algebra $A$ \emph{projectively 
simple} if $\dim_k A = \infty$ and 
every nonzero graded ideal $I$ of $A$ satisfies $\dim_k A/I < \infty$.
Theorem~\ref{xxthm1.1} 
also plays a key role in the study of projectively simple rings in the paper \cite{RRZ}, 
where we prove the following classification result. 
\begin{theorem}\cite[Theorem 0.4 plus Proposition 9.2]{RRZ}
\label{xxthm1.2}
Suppose that $A$ is a connected graded domain which is projectively 
simple, strongly noetherian, is generated in degree 1, and has at least one $k$-point module. 
\begin{enumerate}
\item The canonical map 
$\varphi: A \to B(X, \mc{L}, \sigma)$ of Theorem~\ref{xxthm1.1} is an isomorphism in large degree.
\item $\GK A = 3$ if and only if $X$ is an abelian surface and $\sigma$ 
is the translation automorphism $x \mapsto x + b$ 
for some $b\in X$ such that $\{nb\; |\; n\in {\mathbb Z}\}$ is Zariski-dense in $X$.
\end{enumerate}
\end{theorem}
Some other interesting consequences of Theorem \ref{xxthm1.1}
may be found in Corollary \ref{xxcor4.6} below.

Theorem~\ref{xxthm1.1}(1) follows quite formally from the
Hilbert scheme theorem \cite{AZ2} and the constructions in \cite{ATV}.
The bulk of the work of this article is directed towards the proof
of Theorem~\ref{xxthm1.1}(2).  For this we will develop some results
on the growth of graded rings which are of independent interest.
See Section~\ref{sect2} below for the relevant definitions of growth and
GK-dimension.

Let $R$ be a commutative algebra with algebra automorphism $\sigma$,
and let $Q$ be the skew Laurent polynomial ring $R[t^{\pm 1}; \sigma]$.
Consider a locally finite $\mb{N}$-graded subalgebra
$A = \bigoplus_{n = 0}^{\infty} V_n t^n \subset Q$. We call $A$ a
\emph{big} subalgebra of $Q$ if there exists $n \geq 1$ and a unit
$u \in V_n$ such that $R$ is a localization of its subalgebra
$k[V_n u^{-1}]$ at some multiplicative system of nonzerodivisors.
The key observation needed for the proof of Theorem~\ref{xxthm1.1}(2)
is that all big subalgebras of $Q$ have closely related growth.  
For this result the base field $k$ can be arbitrary.
\begin{proposition}
\label{xxprop1.3}
Suppose $A$ is a big subalgebra of $Q$. If $A$ has finite
GK-dimension (or subexponential growth), then every locally
finite, finitely generated, $\mb{N}$-graded subalgebra
of $Q$ has finite GK-dimension (respectively subexponential growth).
\end{proposition}
The proposition is used in the proof of Theorem~\ref{xxthm1.1}(2) to 
show that the sheaf $\mc{L}$ appearing there is $\sigma$-ample, as follows. 
We shall see that 
$B = B(X, \mc{L}, \sigma)$ and $\varphi(A)$ are 
both big subalgebras of the same skew Laurent polynomial ring $Q$, and 
$\varphi(A)$, being noetherian, has subexponential growth.
Then $B$ has subexponential growth, whereas if $\mc{L}$ were not 
$\sigma$-ample then $B$ would have exponential growth.

Although Proposition~\ref{xxprop1.3} was developed specifically
for its application to the proof of Theorem~\ref{xxthm1.1}(2),
it suggests some further interesting questions unrelated
to Theorem~\ref{xxthm1.1}.

\begin{question}
\label{xxque1.4}
Suppose that $R$ is a commutative algebra with automorphism
$\sigma$ such that there is a locally finite graded big
subalgebra $A \subset R[t^{\pm 1};\sigma]$ with $\GK A < \infty$.
\begin{enumerate}
\item
Does every big locally finite $\mb{N}$-graded subalgebra of
$R[t^{\pm 1};\sigma]$ have the same GK-dimension as $A$?
\item
Must $R[t^{\pm 1};\sigma]$ itself have finite GK-dimension as a
$k$-algebra?
\end{enumerate}
\end{question}

At the end of the paper we will give a positive answer to this
question in a special case.  Specifically, we study the case where
$R[t^{\pm 1};\sigma]$ is the graded ring of fractions of some
twisted ring $B(X, \mc{L}, \sigma)$ where $X$ is integral and
$\mc{L}$ is $\sigma$-ample.

\begin{theorem}
\label{xxthm1.5} Let $X$ be an integral projective scheme
with automorphism $\sigma$, inducing an automorphism of the
function field $K = k(X)$ we also call $\sigma$. Assume that $\mc{L}$
is a $\sigma$-ample sheaf on $X$, and put $d = \GK B(X, \mc{L},
\sigma) < \infty$. Let $Q = K[t^{\pm 1}; \sigma]$.
\begin{enumerate}
\item
Every big locally finite $\mb{N}$-graded subalgebra $A$ of $Q$
has $\GK A = d$.
\item
$d\leq \GK Q \leq d + \dim X$.
\end{enumerate}
\end{theorem}

In general Question \ref{xxque1.4} is still open. As an
application of Theorem \ref{xxthm1.5}, the GK-transcendence
degree of $Q$ is computed in Corollary \ref{xxcor5.8}.

\section{GK-type}
\label{sect2}

Throughout $k$ is a commutative base field, and everything is over
$k$. In particular, an algebra or a ring means an algebra
over $k$. We refer to \cite{KL} for the topic of growth of algebras.

For an algebra $A$, we define the \emph{Gelfand-Kirillov
dimension} (or \emph{GK-dimension}) of $A$ to be
\[
\GK(A) = \sup_V \overline{\lim_{n\to\infty}} \log_n(\dim_k V^n)
\]
where the supremum is taken over all finite-dimensional
$k$-subspaces $V$ of $A$ \cite{GK}.  If $A$ is a finitely
generated algebra then $\GK(A) = \overline{\lim} \log_n(\dim_k V^n)$
for any subspace $V$ of $A$ such that $1 \in V$ and $V$
generates $A$ as a $k$-algebra.  If $A$ is infinitely generated,
then $\GK(A)$ is the supremum of $\GK(B)$ as $B$ ranges over
the finitely generated subalgebras of $A$.

We say that $A$ has \emph{exponential growth} if
\[
\underset{V}{\sup}\ \overline{\lim_{n\to\infty}}
(\dim_k V^n)^{\frac{1}{n}} > 1
\]
where again $V$ ranges over all finite-dimensional
$k$-subspaces of $A$. If $A$ is finitely generated and has
exponential growth then $\overline{\lim} (\dim_k
V^n)^{\frac{1}{n}} > 1$ holds for any subspace $V$ of $A$
with $1 \in V$ which generates $A$ as an algebra.
If $A$ does not have exponential growth then we say it
has \emph{subexponential growth}.  Note that 
if $\GK A < \infty$, then $A$ has subexponential growth.

The algebras of interest in this paper will be $\mb{Z}$-graded
algebras $A = \bigoplus_{n = -\infty}^{\infty} A_n$.  The
$\mb{Z}$-graded algebra $A$ is \emph{$\mb{N}$-graded}
if $A_n = 0$ for all $n < 0$, and \emph{locally finite} if
$\dim_k A_n < \infty$ for all $n \in \mb{Z}$. We call $A$
\emph{finitely graded} if (i) $A$ is locally finite; (ii)
$A$ is $\mb{N}$-graded; and (iii) $A$ is finitely generated
as an algebra over $k$. If $A_0 = k$, $A$ is called \emph{connected}.
If $A$ is locally finite, finitely generated, and $\mb{Z}$-graded,
then
$$\GK A=\overline{\lim_{n\to\infty}} \log_n(\sum_{i=-n}^n\dim_k A_i).$$
If both $A_{\geq 0}$ and $A_{\leq 0}$ are locally finite and
finitely generated, then $\GK A=\max\{\GK A_{\geq 0},\GK A_{\leq 0}\}$.
Given a $\mb{Z}$-graded algebra $A$, for each $d \geq 1$ we have 
the \emph{$d$th Veronese ring} $A^{(d)} = \bigoplus_{n = -\infty}^{\infty} A_{nd}$, 
which is again $\mb{Z}$-graded.

If the set $S$ of all homogeneous nonzerodivisors of
$A$ is an Ore set, then we call the localization $AS^{-1}$
the \emph{graded ring of fractions} of $A$. If the graded
ring of fractions for $A$ exists, and in addition some
element $t$ in $A_1$ is a nonzerodivisor, then the ring of fractions
can be written in the form of a skew Laurent polynomial
ring $Q = R[t^{\pm 1}; \sigma]$ for some ring $R$ with
automorphism $\sigma$.

Now suppose we are given some skew Laurent polynomial
ring $Q = R[t^{\pm 1}; \sigma]$, which we always assume is
$\mb{Z}$-graded with $R = Q_0$ and $t \in Q_1$.  We
want to compare the growth of the various graded
subrings of $Q$.

\begin{definition}
\label{xxdef2.1}
Let $R$ be an algebra, and $\sigma: R \to R$ an algebra
 automorphism.  Let $Q = R[t^{\pm 1}; \sigma]$. Then we say
that $Q$ is of \emph{finite GK-type} if every finitely graded
subalgebra $A \subset Q$ has finite GK-dimension. Similarly,
we say that $Q$ is of \emph{subexponential GK-type} if every
such $A \subset Q$ has subexponential growth.
\end{definition}

We are especially interested in the GK-type of a ring
$Q = R[t^{\pm 1}; \sigma]$ where $R$ is commutative.
In studying the GK-type we may restrict our attention
to those graded subalgebras $A$ of $Q$ which are
largest in some sense; intuitively, this means those
$A$ for which the graded ring of fractions of $A$ is
$Q$ and not some smaller ring.  Since graded rings of
fractions do not always exist, we will use instead
the following notion of largeness for subrings of $Q$,
which we already mentioned in the introduction.

\begin{definition}
\label{xxdef2.2}
Let $Q = R[t^{\pm 1}; \sigma]$ where $R$ is a commutative
algebra with algebra automorphism $\sigma$.
A locally finite ($\mb{N}$-)graded subalgebra $A =
\bigoplus_{n = 0}^{\infty} V_n t^n \subset Q$ is called
a \emph{big} subalgebra of $Q$ if there is some $n \geq 1$
and an element $u \in V_n$ such that (i) $u$ is a unit of
$R$; and (ii) $R$ is a localization of its subalgebra
$k[V_n u^{-1}]$ at some multiplicative system of
nonzerodivisors of $k[V_n u^{-1}]$.
\end{definition}

In special cases, for example if $R$ is a field, the
concept of a big subalgebra can be formulated in a more
intuitive way.  The proof of the following alternative
characterizations is straightforward.

\begin{lemma}
\label{xxlem2.3}
Let $Q = K[t^{\pm 1}; \sigma]$ with $K$ a field over $k$.
Let $A$ be a locally finite graded subalgebra of $Q$.
\begin{enumerate}
\item
$A$ is big in $Q$ if and only if the smallest graded
simple subring of $Q$ containing $A$ is
$K[(t^s)^{\pm 1}; \sigma^s]$ for some $s \geq 1$.
\item
If the set of homogeneous nonzerodivisors of $A$ is
an Ore set, then $A$ is big in $Q$ if and only if
the graded ring of fractions for $A$ is
$K[(t^s)^{\pm 1}; \sigma^s]$ for some $s \geq 1$.
\end{enumerate}
\end{lemma}

Our main result concerning GK-type is
Proposition~\ref{xxprop1.3}, which shows that if $R$ is
commutative, then to prove that $Q = R[t^{\pm 1}; \sigma]$
is of finite (or subexponential) GK-type, it suffices to
find a single big subalgebra of $Q$ of finite GK-dimension
(respectively subexponential growth). In the proof below,
for any subset $W \subset R$ and automorphism $\tau$ of
$R$ we use the notation $W^{\tau} = \tau(W)$.

\begin{proof}[Proof of Proposition~\ref{xxprop1.3}]
Let $A = \bigoplus_{n = 0}^{\infty} V_n t^n$ be a big
subalgebra of $Q$, and let $B = \bigoplus_{n = 0}^{\infty}
W_n t^n$ be an arbitrary finitely graded subring of $Q$. We
want to show that
if $A$ has finite GK-dimension or subexponential growth,
then the same is true of $B$.  
If $B$ is generated as an algebra in degrees $\leq m$, then 
the ring $k + B_{\geq 1}$ is also generated in degrees $\leq m$, 
and has the same growth as $B$.  Replacing $B$ by 
this algebra we way assume that $B$ is connected.  Now 
$B$ is contained in the algebra generated by $k + Wt$ where
$$W = k + \sum_{n=1}^m W_n.$$
Replacing $B$ by this larger algebra we may assume that
$B_n= (B_1)^n$ for $n\geq 1$, so that $B$ is generated in degree $1$.

By the hypothesis that $A$ is big we may find $n \geq 1$
and a unit $u \in V_n \subset R$ such that $R$ is a
localization of $k[V_n u^{-1}]$.  Consider the Veronese
subrings $A^{(n)}, B^{(n)}$, which are graded
subrings of $Q^{(n)} \cong R[(t^n)^{\pm 1}; \sigma^n]$. Since
$B$ is finitely generated over $B^{(n)}$, $B^{(n)}$ and
$B$ have the same growth; since $A^{(n)} \subset A$ as
ungraded rings, $A^{(n)}$ grows no faster than $A$.
Also $A^{(n)}$ is big in $Q^{(n)}$.  Thus we may replace
$A$ and $B$ by their $n$th Veronese subrings without harm.
Then setting $V = V_1$, $W = W_1$, and $U = V u^{-1}$, we
have $1 \in U$ and $R$ is a localization of $k[U]$.

Thinking of the elements in a basis of $W$ as fractions in
$k[U]$ and putting them over a common denominator, we see
that we can find an integer $f$ and a unit $u' \in R$ such
that
\[
 W \subset u' U^f =  z V^f,
\]
where $z = u' u^{-f}$ is also a unit.  Hence we have
\[
W W^{\sigma} \cdots W^{\sigma^n}\subset z\cdots z^{\sigma^n}
(V\cdots V^{\sigma^n})^f
\]
and thus
\[
\dim_k W\cdots W^{\sigma^n} \leq  \dim_k (V\cdots V^{\sigma^n})^f
\leq  (\dim_k V\cdots V^{\sigma^n})^f.
\]
Because $B_{n+1}=(W\cdots W^{\sigma^n})t^{n+1}$ and
$A_{n+1} \supset (V\cdots V^{\sigma^n})t^{n+1}$, it
follows that if $A$ has finite GK-dimension (or
subexponential growth), then $B$ also has finite GK-dimension
(respectively subexponential growth), as we wished to show.
\end{proof}

The following consequence is immediate.

\begin{corollary}
\label{xxcor2.4} Let $Q=R[t^{\pm 1};\sigma]$ where $R$ is a
commutative algebra. Let $A$ be a finitely graded
subalgebra of $Q$. If $A$ has exponential growth, then every
finitely graded big subalgebra of $Q$ has exponential growth.
\end{corollary}

If $R$ is not commutative but only finite over its center,
we can understand the GK-type of $Q = R[t^{\pm 1}; \sigma]$
by reducing to the commutative case, using the following
result.  We include this result for the interested reader,
but we will not need to use it later on.

\begin{proposition}
\label{xxprop2.5}
Let $R$ be an algebra with automorphism $\sigma$. Let
$S$ be a subalgebra of $R$ such that $S$ is in the
center of $R$, $\sigma$ restricts to an automorphism of
$S$, and $R$ is a finite $S$-module.  Then $Q =
R[t^{\pm 1}; \sigma]$ is of finite (or subexponential)
GK-type if and only if $Q' = S[t^{\pm 1}; \sigma]$ is of
finite (respectively subexponential) GK-type.
\end{proposition}

\begin{proof}
If $Q$ is of finite or subexponential GK-type, then the
same obviously holds for the subring $Q'$ of $Q$.

For the converse, let $B$ be an arbitrary finitely
graded subalgebra of $Q$.  We will show that there is a
finitely graded subalgebra $A$ of $Q'$ with
$\dim_k B_n < d \dim_k A_n$ for some constant $d > 0$
and all $n \geq 0$.  This is sufficient to complete
the proof.

By the same reductions as in the proof of Proposition~\ref{xxprop1.3}, 
we may assume that $B$ is connected and generated in degree one, say by 
$Wt$.  For the rest of this proof, for any finite-dimensional
$k$-subspace $Z$ of $R$ we write $Z_0 = k$ and
$Z_n = Z Z^{\sigma} \cdots Z^{\sigma^{n-1}}$ for all $n \geq 1$.
With this notation we have $B = \bigoplus_{n \geq 0} W_n t^n$.

Write $R = \sum_{i =1}^d S r_i$ for some fixed elements
$r_i \in R$. Let $T$ be the finite-dimensional $k$-span
of $\{r_i\}$. We may choose finite-dimensional
subspaces $V,U\subset S$ with $1 \in U$ such that $W\subset VT$
and $TT^\sigma\subset UT$.
Now one proves by induction
that
\[
T T^{\sigma} \cdots T^{\sigma^{n-1}} \subseteq T U U^{\sigma}
\cdots U^{\sigma^{n-2}}
\subseteq T U U^{\sigma} \cdots U^{\sigma^{n-1}}
\]
for all $n \geq 1$.  In other words, $T_n \subset TU_n$
for all $n \geq 1$.  Then we have
$$W_n \subseteq T_n V_n \subseteq T U_n V_n.$$
Hence $\dim_k W_n \leq d \dim_k U_n V_n$ for all
$n \geq 0$, where $d = \dim_k T$. Define $A \subset Q'$
to be the connected graded ring generated in degree $1$ by
$UVt$. Then $A$ is the desired subring of $Q'$ with
$\dim_k B_n < d \dim_k A_n$ for all $n \geq 0$.
\end{proof}

The similarity in the growth of graded subrings of
$Q = R[t^{\pm 1}; \sigma]$ which is exhibited by
Propositions~\ref{xxprop1.3} and \ref{xxprop2.5}
does not necessarily extend to the case where $R$ is not finite
over its center, as the following example shows.

\begin{example}
\label{xxex2.6}
Let $A$ be the graded algebra defined by generators
and relations as $A = k \langle x,y,z \rangle/
(xy-yx-z^2, xz-zx, yz-zy)$.  Then $A$ is one of the
standard examples of a $3$-dimensional Artin-Schelter
regular algebra, and $\GK A = 3$. The graded ring of
fractions of $A$ is $Q = D[z^{\pm 1};id_D]$ where $D$ is the
first Weyl skew-field.
The division ring $D$ contains two elements $f, g \in D$
which generate a free algebra \cite[Theorem 8.17]{KL}.
Thus the subring $A' = k [Vz] \subset Q$, where $V$
is the $k$-span of $\{1,f,g\}$, has exponential growth,
even though $Q$ is the graded ring of fractions of the
subring $A$ of finite GK-dimension.
\end{example}

\section{Twisted rings over non-$\sigma$-ample $\mc{L}$}
\label{sect3}

Throughout this section \emph{let $k$ be an algebraically closed field, 
$X$ be a projective $k$-scheme, $\sigma$ an automorphism of $X$, and 
$\mc{L}$ an invertible
sheaf over $X$.} We will show some preliminary results about
twisted homogeneous coordinate rings (which we call  
twisted rings for brevity) for non-$\sigma$-ample
invertible sheaves $\mc{L}$.  

First we recall the definitions from \cite{AV}. For any sheaf
$\cal{F}$ on $X$, we use the notation $\cal{F}^{\sigma}$ for the
pullback $\sigma^{*}(\cal{F})$. Set $\mc{L}_0 = \mc{O}_X$ and
$\cal{L}_n = \cal{L} \otimes \cal{L}^{\sigma} \otimes \dots
\otimes \cal{L}^{\sigma^{n-1}}$ for all $n \geq 1$.
The \emph{twisted ring} associated to the triple $(X,\mc{L},
\sigma)$ is
$$B(X, \mc{L}, \sigma)=\bigoplus_{n = 0}^{\infty}
\HB^0(X, \mc{L}_n),$$
which has a natural $\mb{N}$-graded algebra structure given
by the multiplication rule
$$f \cdot g = f \otimes (\sigma^m)^*(g)$$
for $f \in \HB^0(X, \mc{L}_m), g \in \HB^0(X, \mc{L}_n)$. The
invertible sheaf $\mc{L}$ is called \emph{$\sigma$-ample} if
for any sheaf $\cal{F} \in \coh X$, $\HB^q(\cal{F} \otimes
\cal{L}_n) = 0$ for all $q > 0$ and $n \gg 0$. Other details
can be found in \cite{AV,Ke1}.

Keeler showed that $\sigma$-ampleness depends on a property
of the automorphism $\sigma$ called quasi-unipotence \cite{Ke1}.
To define this, let $\Num X$ be the quotient of the Picard group
of $X$ by the subgroup of divisors numerically equivalent to $0$;
then $\sigma$ induces an action of the finitely generated free
abelian group $\Num X$, and $\sigma$ is called
\emph{quasi-unipotent} if the eigenvalues of this action are
roots of unity. We will not need the explicit definition of
quasi-unipotence below, but only the criterion that Keeler
proved:  if $\sigma$ is quasi-unipotent, then an
invertible sheaf $\mc{L}$ is $\sigma$-ample if and only if
$\mc{L}_n$ is ample for some $n$, while if $\sigma$ is not
quasi-unipotent then no invertible sheaf on $X$ is
$\sigma$-ample \cite[Theorem 1.3]{Ke1}.  When $\mc{L}$ is
$\sigma$-ample then $B(X, \mc{L}, \sigma)$ is a noetherian
ring with well-understood properties.  Comparatively little
is known about twisted rings $B(X, \mc{L}, \sigma)$ where
$\sigma$ is not quasi-unipotent. The existing knowledge
about this case is again due to Keeler, as follows.

\begin{lemma} \cite[Proposition 3.7]{Ke1}
\label{xxlem3.1}
Suppose $\sigma$ is non-quasi-unipotent. Let $\mc{L}$ be an
ample invertible sheaf. Then there is an integer $m_0$ with the property
that for any given $m \geq m_0$, one can find constants $C > 0$ and 
$r > 1$ such that the graded pieces of $B: =
B(X, \mc{L}^{\otimes m}, \sigma)$ satisfy $\dim_k B_n \geq
Cr^n$. In particular, it follows
that for any $m \geq m_0$ the ring $B$ is not noetherian.
\end{lemma}

\begin{remark}
\label{xxerm3.2}
Keeler defines a locally finite $\mb{N}$-graded ring $A =
\bigoplus_{n = 0}^{\infty} A_n$ to have exponential growth
if $\overline{\lim} \big(\sum_{i \leq n}
\dim_k A_i\big)^{1/n} > 1$ \cite[Equation (3.3)]{Ke1}.  This
definition agrees with ours (see Section~\ref{sect2}) only if
$A$ is finitely generated.  An infinitely generated algebra
$A$ could well have graded pieces $A_n$ growing exponentially
in size with $n$, yet have finite GK-dimension under our
definitions if it is a directed union of finitely generated
subrings of the same finite GK-dimension.  Keeler asserts in
\cite[Proposition 3.7]{Ke1} that if $\sigma$ is not
quasi-unipotent and $\mc{L}$ is ample, then for $m \geq m_0$
the ring $B = B(X, \mc{L}^{\otimes m}, \sigma)$ has exponential
growth. What is actually proven is the exponential growth of
graded pieces $\dim_k B_n \geq Cr^n$ for some $r > 1$ and
$C > 0$; we take this as the conclusion of Keeler's
proposition in Lemma~\ref{xxlem3.1} above, since it is not
enough to imply that $B$ has exponential growth under our
definition. To conclude that $B$ is not noetherian for
$m \geq m_0$, Keeler uses \cite[Theorem 0.1]{SZ} which
states that a locally finite graded ring of exponential growth
cannot be noetherian.  This conclusion is unaffected by
the different definitions of exponential growth, since a
noetherian locally finite $\mb{N}$-graded algebra would have to
be finitely generated as an algebra in any case.
\end{remark}

We will generalize Lemma~\ref{xxlem3.1} to show that if
$\sigma$ is not quasi-unipotent but some $\mc{L}_n$ is
\emph{very} ample, then the ring $B(X, \mc{L}, \sigma)$
already has exponential growth.  As we see from the preceding
remark, it will be useful first to study finite generation of
twisted rings.  For this, the theory of Castelnuovo-Mumford
regularity will provide a convenient tool. 
A coherent sheaf $\mc{F}$ on $X$ is said to be \emph{$m$-regular}
(with respect to $\mc{L}$) if $\HB^i(\mc{F} \otimes
\mc{L}^{\otimes m-i}) = 0$ for all $i \geq 1$.  It is a fact
that if a sheaf $\mc{F}$ is $m$-regular, then it is also
$m'$-regular for all $m' \geq m$.  The \emph{regularity} of
$\mc{F}$ is the smallest integer $m$ for which $\mc{F}$ is
$m$-regular, and is denoted by $\reg \mc{F}$ once the sheaf
$\mc{L}$ is fixed.

Now we will see that when $\mc{L}$ is ``ample enough'' the
ring $B(X, \mc{L}, \sigma)$ is a finitely generated
algebra, regardless of whether $\mc{L}$ is $\sigma$-ample
or not.

\begin{proposition}
\label{xxprop3.3}
Let $\mc{L}$ be an ample invertible sheaf on $X$.  Consider
the ring $B = B(X, \mc{L}^{\otimes m}, \sigma)$.  Then 
there is some $m_0 \geq 1$ such that if $m \geq m_0$, 
then the multiplication map $B_1 \otimes B_{n-1} \to B_{n}$ 
is surjective for all $n \geq 1$.  In other
words, for $m \geq m_0$ 
the ring $B$ is generated in degrees $0$ and $1$.
\end{proposition}

\begin{proof}
Fix some arbitrary very ample sheaf $\mc{P}$ on $X$ and
measure Castelnuovo-Mumford regularity with respect to
$\mc{P}$ for the rest of the proof, which we break up
into several steps.

\emph{Step 1}.
There exists a constant $C_1$ such that for all coherent
sheaves $\mc{F}, \mc{M}$ on $X$ with $\mc{M}$ invertible,
we have
\[
\reg (\mc{F} \otimes_{\mc{O}_X} \mc{M}) \leq \reg \mc{F}
+ \reg{\mc{M}} + C_1.
\]
This is proved by Keeler in \cite{Ke2}.

\emph{Step 2}.
There is a constant $C_2$ such that $\reg \mc{N} \leq C_2$
for all ample invertible sheaves $\mc{N}$ on $X$.  This is
a consequence of a vanishing Lemma of Fujita
\cite[p. 520, Theorem 1]{Fj}.

\emph{Step 3}.
Let $\mc{M}$ be an ample invertible sheaf on $X$.  Then
for any fixed $j \in \mb{Z}$ we have
$\HB^i(\mc{P}^{\otimes j} \otimes \mc{M}^{\otimes n}) = 0$
for $n \gg 0$ and $i > 0$ by Serre vanishing.  It follows
that $\lim_{n \to \infty} \reg \mc{M}^{\otimes n}
= - \infty$.

\emph{Step 4}.
For each $m \geq 1$, since $\mc{L}^{\otimes m}$ is
generated by its global sections we may consider the
exact sequence
\[
0 \to \mc{K}_m \to \HB^0(\mc{L}^{\otimes m}) \otimes
\mc{O}_X \to \mc{L}^{\otimes m} \to 0.
\]
Then there is a constant $C_3$ such that $\reg \mc{K}_m
\leq C_3$ for all $m \geq 1$.  To see this, note that for
$m \gg 0$, $\mc{L}^{\otimes m}$ is $0$-regular, by step 3.
Then for $m \gg 0$ \cite[Lemma 3.1]{Ar} may be applied
to show that $\reg \mc{K}_m \leq \reg \mc{O}_X + 1$.

\emph{Step 5}.
Consider for all $m \geq 1$, $n \geq 1$ the sheaf
$\mc{G}_{m,n} = \mc{K}_m \otimes
(\mc{L}^{\otimes m}_{n-1})^{\sigma}$, where $\mc{K}_m$
is defined as in step 4. Since $\mc{G}_{m,n} = \mc{K}_m
\otimes (\mc{L}^{\sigma})^{\otimes m} \otimes \mc{N}$ where
$\mc{N} = (\mc{L}^{\sigma^2} \otimes \cdots \otimes
\mc{L}^{\sigma^n})^{\otimes m}$ is ample, we see from
steps 1,2, and 4 that
$$\reg \mc{G}_{m,n} \leq 2 C_1 +
C_2 + C_3 + \reg (\mc{L}^{\sigma})^{\otimes m}.$$
Then by step 3 applied with $\mc{M} = \mc{L}^{\sigma}$,
we may choose $m_0$ such that $\reg \mc{G}_{m,n} \leq 1$
for all $m \geq m_0$ and all $n \geq 1$. In particular,
we have $\HB^1(\mc{G}_{m,n}) = 0$ for $m \geq m_0$ and
all $n \geq 1$. Now tensoring the exact sequence of
step 4 with $(\mc{L}^{\otimes m}_{n-1})^{\sigma}$, we
get an exact sequence
\[
0 \to \mc{G}_{m,n} \to \HB^0(\mc{L}^{\otimes m})
\otimes (\mc{L}^{\otimes m}_{n-1})^{\sigma} \to
\mc{L}^{\otimes m}_n \to 0.
\]
Taking the long exact sequence in cohomology gives
\[
0 \to \HB^0(\mc{G}_{m,n}) \to \HB^0(\mc{L}^{\otimes m})
\otimes \HB^0((\mc{L}^{\otimes m}_{n-1})^{\sigma}) \to
\HB^0(\mc{L}^{\otimes m}_n) \to \HB^1(\mc{G}_{m,n}) \to \dots
\]
and thus the map $\HB^0(\mc{L}^{\otimes m}) \otimes
\HB^0((\mc{L}^{\otimes m}_{n-1})^{\sigma}) \to
\HB^0(\mc{L}^{\otimes m}_n)$ is a surjection for
$m \geq m_0$ and  $n \geq 1$.  This is exactly the desired
claim that the multiplication map $B_1 \otimes B_{n-1} \to
B_n$ is surjective for all $n \geq 1$, for $m \gg 0$.
\end{proof}

Next we will prove a needed technical lemma, which is
standard in case $X$ is integral, but perhaps not so
well-known for general projective schemes $X$. Given any
projective scheme $X$, let $\mc{K}$ be the \emph{sheaf of
total quotient rings} of $X$ \cite[p.~141]{Ha} and let
$R$ be $\HB^0(\mc{K})$, which is called the
\emph{function ring} of $X$ following Nakai. If $X$
is integral, then $R$ is the function field $k(X)$
and $\mc{K}$ is simply the constant sheaf on the function
field $k(X)$. If $X$ is not integral, $\mc{K}$ may not be
a constant sheaf. Nakai proved that every invertible sheaf
$\mc{L}$ on $X$ is equal to a subsheaf of $\mc{K}$
\cite[Theorem~4]{Na}. The automorphism $\sigma$ induces
an automorphism of $R$, which we also call $\sigma$.
If we fix an embedding $\mc{L} \subset \mc{K}$ (following
\cite[Theorem~4]{Na}), then for each $i \geq 0$ using the natural 
isomorphism $\mc{K}^{\sigma^i} \cong \mc{K}$ we get an embedding
$\mc{L}^{\sigma^i} \subset \mc{K}$.  Also, given any two invertible 
subsheaves $\mc{M}, \mc{N} \subset \mc{K}$, by thinking of 
them as Cartier divisors it is clear that $\mc{M} 
\otimes \mc{N} \subset \mc{K}$ \cite[II.6.13]{Ha}.  Thus we have an embedding 
$\mc{L}_n \subset \mc{K}$ for each $n \geq 0$ and 
thus an inclusion $\HB^0(\mc{L}_n) \subset R$.  Finally, 
by adding a placeholder variable $t$,
we may think of the twisted ring $B = B(X, \mc{L}, \sigma)$
as the explicit subring
$\bigoplus_{n = 0}^{\infty} \HB^0(\mc{L}_n)t^n$ of
$Q = R[t^{\pm 1}; \sigma]$.

\begin{lemma}
\label{xxlem3.4}
Let $\mc{L} \subset \mc{K}$ be a very ample invertible sheaf,
and let $V \subset \HB^0(\mc{L}) \subset R$ be a
finite-dimensional $k$-subspace such that the global sections in
$V$ generate the sheaf $\mc{L}$.  Suppose the corresponding
map to projective space $\phi: X \to \mb{P}(V)$ is a closed
embedding. Then there exists a unit $u \in V$ such that $R$
is the total quotient ring of $k[u^{-1}V]$.
\end{lemma}

\begin{proof}
Let $d +1= \dim_k V$, and write $\mb{P}^d = \mb{P}(V)$.  Let
$\mc{I}$ be the ideal sheaf on $\mb{P}^d$ defining the closed
subscheme $X$. Since $\mc{L} \cong \phi^*(\mc{O}(1))$, we
have a homomorphism of graded rings
\[
k[z_0, \dots, z_d] = \textstyle \bigoplus_{n = 0}^{\infty}
\HB^0(\mb{P}^d, \mc{O}(n)) \overset{\theta}{\lra} \textstyle
\bigoplus_{n = 0}^{\infty} \HB^0(X, \mc{L}^{\otimes n})t^n
\subset R[t^{\pm 1}].
\]
The kernel of $\theta$ is $I = \bigoplus_{n = 0}^{\infty}
\HB^0(\mb{P}^d, \mc{I} \otimes \mc{O}(n))$, and the image
is the ring $A = \bigoplus_{n = 0}^{\infty} V^n t^n$.
Then $A \cong k[z_0, \dots, z_d]/I$, and $A$ has no
finite-dimensional ideals \cite[Exercise II.5.10]{Ha}.  Thus
$A_1$ contains some nonzerodivisor $z = ut$.
Let $A_{(z)}$ be the degree $0$ piece of the localization of
$A$ at $z$.  By  \cite[p. 298]{Na}, we see that $R$ is the
total quotient ring of $A_{(z)}$.  Thus $R$ is the total
quotient ring of $k[Vu^{-1}] \subset R$, where $u \in V$
is a unit in $R$, as required.
\end{proof}

We are now ready to prove our generalization of
Lemma~\ref{xxlem3.1}.

\begin{proposition}
\label{xxprop3.5}
Let $\sigma$ be a non-quasi-unipotent automorphism of $X$.
\begin{enumerate}
\item Every big locally finite $\mb{N}$-graded subalgebra of
$R[t^{\pm 1}; \sigma]$ has exponential growth.
\item Let $\mc{L}$ be an invertible sheaf on $X$ such
that $\mc{L}_n$ is very ample for some $n \geq 1$.
Then $B(X, \mc{L}, \sigma)$ has exponential growth and
is not noetherian.
\end{enumerate}
\end{proposition}
\begin{proof}
(1) Let $\mc{L}$ be a very ample invertible sheaf on
$X$ with an embedding $\mc{L} \subset \mc{K}$.
By Lemma~\ref{xxlem3.1}, for $m \gg 0$ the ring
$B = B(X, \mc{L}^{\otimes m}, \sigma)$ has the property
that $\overline{\lim} (\dim_k B_n)^{1/n} > 1$.  By
Proposition~\ref{xxprop3.3} above, for $m \gg 0$ we also
know that $B$ is finitely generated; together these facts
imply that $B$ has exponential growth for $m \gg 0$.
Setting $\mc{L}' = \mc{L}^{\otimes m} \subset \mc{K}$, we
may write $B = \bigoplus_{n = 0}^{\infty}
\HB^0(X, \mc{L}'_n) t^n \subset Q:=R[t^{\pm 1}; \sigma]$.
By Corollary \ref{xxcor2.4}, every big finitely graded subalgebra of $Q$
has exponential growth.  Since every big locally finite 
$\mb{N}$-graded subalgebra of $Q$ 
contains a finitely graded big subalgebra, the result follows.

(2) Now let $\mc{L}$ be an invertible sheaf on $X$ such
that $\mc{L}_n$ is very ample for some $n$. Fix an embedding
$\mc{L} \subset \mc{K}$ and let $B = B(X, \mc{L}, \sigma) \cong
\bigoplus_{n = 0}^{\infty} \HB^0(\mc{L}_n) t^n \subset Q$.
Setting $V_n = \HB^0(\mc{L}_n) \subset R$, since $\mc{L}_n$ is very ample 
the sections $V_n$ generate $\mc{L}_n$, and also the induced map to 
projective space $X \to \mb{P}(V_n)$ is a closed embedding.
Then it follows by Lemma~\ref{xxlem3.4} that there is a unit $u \in V_n$ such
that $R$ is the total quotient ring of $k[u^{-1}V_n]$.  So 
by definition $B$ is a big subalgebra of $Q$.  
Then $B$ has exponential growth by part (1). By \cite[Theorem~0.1]{SZ},
$B$ is not noetherian.
\end{proof}

We note that Proposition~\ref{xxprop3.5}(2) gives a partial
answer to \cite[Question 3.8]{Ke1}, which asks if given a
non-quasi-unipotent automorphism $\sigma$ and an ample
invertible sheaf $\mc{L}$, whether $B(X, \mc{L}, \sigma)$ must
have exponential growth.  We do not know in the non-quasi-unipotent
case if given an ample $\mc{L}$ there is always some $n \geq 1$
such that $\mc{L}_n$ is very ample.

\section{The canonical map $\varphi$}
\label{sect4}

In this section we assume that $k$ is algebraically closed since
we will quote some results from \cite{ATV} and \cite{AZ2} which
are under this assumption. A module means a right module.
We begin this section with a review of the results of
\cite[\S 3]{ATV}.  We refer the reader to
that paper for more complete details.

Let $A$ be a graded algebra generated in degree $1$ and let
$\dim_k A_1 = d+1$. Let $T=k \langle x_0, x_1, \dots, x_d \rangle$
be the free algebra, and fix a presentation $A \cong T/I$.
Let $V = (T_1)^* = \Hom_k(T_1, k)$.  Any $f \in T_n$ determines an
element in $(V^{\otimes n})^*$, and thus an element $\wt{f} \in
\HB^0(\mb{P}^{\times n}, \mc{O}(1,1,\dots, 1))$, where
$\mb{P}^{\times n}$ is an abbreviation for $(\mb{P}(V))^{\times n}$.
Then $\wt{f}$ determines a vanishing locus in $\mb{P}^{\times n}$.
Now for each $n \geq 0$, we let $X_n \subset \mb{P}^{\times n}$ be
the common locus of vanishing of all $\{\wt{f} \mid f \in I_n \}$.
Also for each $n \geq 0$ let $i_n: X_n \to \mb{P}^{\times n}$ be
the inclusion map, and set $\mc{M}_n = i_n^*\mc{O}(1,1, \dots 1)$.
Now for any $1 \leq m < n$ we define a morphism of schemes
$\phi_{n,m}: X_n \to X_m$  by restricting the morphism $\pi_1:
\mb{P}^{\times n} \to \mb{P}^{\times m}$ which is the projection
onto the first $m$ coordinates.  We also define a morphism
of schemes $\psi_{n,m}: X_n \to X_m$ by restricting the morphism
$\pi_2: \mb{P}^{\times n} \to \mb{P}^{\times m}$ which is the
projection onto the last $m$ coordinates.  We use the
abbreviations $\phi_{n+1, n} = \phi_n$ and $\psi_{n+1, n} = \psi_n$.

The schemes just constructed are closely related to the point
modules for the ring $A$. For any commutative algebra $R$ and
$n \geq 0$, a \emph{truncated $R$-point module of length $n+1$}
is a graded cyclic right $A \otimes R$-module $M =
\bigoplus_{i = 0}^n M_i$ where $M_0 = R$ and each $M_i$ for
$1 \leq i \leq n$ is a locally free $R$-module of rank $1$.
As we defined in the introduction, an \emph{$R$-point module}
is defined similarly, by letting $n = \infty$.  A $k$-point
module is simply called a point module.  Then by
\cite[Proposition 3.9]{ATV}, the scheme $X_n$ represents the
functor from commutative $k$-algebras to sets which sends an
algebra $R$ to the set of isomorphism classes of truncated
$R$-point modules of length $n+1$. Moreover, in case the maps
$\phi_n: X_{n+1} \to X_n$ are isomorphisms for all $n \geq n_0$,
then $X = X_{n_0}$ is the point scheme for $A$ as defined in
the introduction; that is, $X$ represents the functor sending
an algebra $R$ to the set of isomorphism classes of $R$-point
modules.

In the next lemma we gather together some further formal
consequences of the definitions above, which are also largely
a restatement of results in \cite{ATV}.

\begin{lemma}
\label{xxlem4.1}
Keep the notations introduced above.
\begin{enumerate}
\item
There is a graded ring homomorphism
\[
\wt{\theta}: T \cong \bigoplus_{n = 0}^{\infty}
\HB^0(\mb{P}^{\times n}, \mc{O}(1,1, \dots, 1))
\to B' = \bigoplus_{n = 0}^{\infty} \HB^0(X_n, \mc{M}_n)
\]
defined by pullback of sections, and $\wt{\theta}$ factors
through $I$ to give a map $\theta: A \to B'$.  An element
$x \in T_n = \HB^0(\mb{P}^{\times n}, \mc{O}(1,1, \dots, 1))$
is in $(\ker \wt{\theta})_n$ if and only if $i_n^*(x) = 0$.
\item
Let $J = \ker \theta \subset A$.  Then for each $n \geq 0$,
$a \in J_n$ if and only if $(M_0) a = 0$ for all truncated
$R$-point modules $M$ of length $n+1$ over $A$, for all
commutative $k$-algebras $R$.
\end{enumerate}
\end{lemma}

\begin{proof}
(1) The ring structure on $B'$ is described in
\cite[p. 48]{ATV}, but we will also describe it briefly
here.  Given $m, n \geq 1$ we have an isomorphism of
sheaves on $X_{n+m}$
\[
\phi_{m+n,m}^*(\mc{M}_m) \otimes \psi_{m+n,n}^*(\mc{M}_n)
\cong \mc{M}_{m+n},
\]
and thus a map
\[
\HB^0(X_m, \mc{M}_m) \otimes \HB^0(X_n, \mc{M}_n)
\overset{\phi^* \otimes
\psi^*}{\lra}
\HB^0(X_{m+n}, \mc{M}_{m+n}).
\]
This gives the multiplication maps $B'_m \otimes B'_n
\to B'_{m+n}$ for each $m, n$.

The map $\wt{\theta}$ is defined in degree $n$ by pulling
back sections via the embedding $i_n: X_n \to \mb{P}^{\times n}$.
That this $\wt{\theta}$ defines a ring homomorphism
factoring through $I$ is \cite[Proposition 3.20]{ATV}.
It is immediate from the definition of $\wt{\theta}$ as a
pullback that $x \in (\ker \wt{\theta})_n$ if and only if
$i_n^* (x) = 0$.

(2) By \cite[Proposition 3.9]{ATV}, for every commutative
algebra $R$ there is a bijection between
$\Hom_{k-\operatorname{schemes}}(\spec R, X_n)$ and the
set of isomorphism classes of truncated $R$-point modules
of length $n+1$ for $A$.  One way to make the correspondence
explicit is to use the universal truncated point module
of length $n+1$ (see \cite[p. 47]{ATV}).  This is a
coherent sheaf on $X_n$ given by
\[
\mc{N} = \bigoplus_{j = 0}^n \mc{N}_j
\quad \text{where}\quad
\mc{N}_j=i_n^* \mc{O}(1,1, \dots, \overset{j}{1}, 0,0, \dots, 0).
\]
Note that $\mc{N}_0 = \mc{O}_{X_n}$.
Fix some $j$ with $1 \leq j \leq n$.  We have maps
\[
A_j \overset{\theta_j}{\to}
\HB^0(X_j, \mc{M}_j) \overset{\phi_{n,j}^*}{\to}
\HB^0(X_n, \mc{N}_j)
\]
and thus we get a morphism of sheaves on $X_n$
\[
\mc{N}_0 \otimes_k A_j \to \mc{O}_{X_n}
\otimes_k \HB^0(X_n, \mc{N}_j) \to \mc{N}_j.
\]
Now given a commutative algebra $R$ and a
$k$-morphism $\rho: \spec R \to X_n$, we can pull
back via $\rho$, and we get a morphism of sheaves on $\spec R$
\begin{equation}
\label{E4.2}
\rho^*(\mc{N}_0) \otimes_k A_j \to \rho^*(\mc{N}_j).
\end{equation}
Then the truncated $R$-point module of length $n+1$
corresponding to the map $\rho$ is simply $M = \rho^* \mc{N}$;
note that $M_0 = R$ and that the maps \eqref{E4.2} determine
the $A$-module structure on $M$.

There is a canonical generator $1 \in R = M_0$ for the
truncated point module $M$.  Then given $a \in A_n$,
if we lift $a$ to some $\wt{a} \in T_n$, then we see
that $(M_0) a = 0$ if and only if $\rho^* i_n^*(\wt{a}) = 0$.

Now we see from part (1) that $a \in (\ker \theta)_n$
if and only if given any lift  $\wt{a} \in T_n$ of $a$,
we have $i_n^*(\wt{a}) = 0$.  But the section $i_n^*(\wt{a})
\in \HB^0(X_n, \mc{M}_n)$ is zero if and only if
$\rho^* i_n^*(\wt{a}) = 0$ for all commutative algebras
$R$ and all maps $\rho: \spec R \to X_n$, if and only if
$(M_0) a = 0$ for all truncated $R$-point modules $M$ of
length $n+1$ by the explicit correspondence given above.
\end{proof}

We need one last ingredient for the proof of the main
theorem.  The following result of Artin and Stafford
shows that under certain conditions, a subring $A$ of a
twisted ring $B$ must be equal to $B$ in large degree.

\begin{lemma}
\cite[Theorem 4.1]{AS}
\label{xxlem4.3}
Suppose that $B = B(X, \mc{L}, \sigma)$ where $\mc{L}$
is $\sigma$-ample.  Let $A$ be a graded subring of $B$,
and suppose that for $n \gg 0$ the sections in $A_n
\subseteq \HB^0(\mc{L}_n)$ generate the sheaf $\mc{L}_n$,
and that the map to projective space $X \to \mb{P}(A_n)$
determined by a $k$-basis for $A_n$ is a closed embedding.
Then $A_n = B_n$ for $n \gg 0$.
\end{lemma}

\begin{proof}
This is a slight restatement of \cite[Theorem 4.1]{AS},
which requires in addition to our hypotheses that $\mc{L}$
is also $\sigma^{-1}$-ample, and that $A$ has subexponential
growth.  But since $\mc{L}$ is $\sigma$-ample, $\mc{L}$ must
also be $\sigma^{-1}$-ample \cite[Corollary 5.1]{Ke1}, and
$B$ has finite GK-dimension \cite[Theorem 6.1]{Ke1}, so $A$
does as well.
\end{proof}

We can now prove our most general result about when maps
to twisted rings exist and are surjective in large degree.

\begin{theorem}
\label{xxthm4.4}
Let $A$ be a connected graded algebra of subexponential
growth which is generated in degree $1$, and let $X_n, i_n, 
\phi_n, \psi_n,$ and $\mc{M}_n$ be defined for each $n \geq 0$
as at the beginning of this section.  Assume that there exists some $n_0$ such that 
for all $n \geq n_0$, both $\phi_n$ and $\psi_n$ are isomorphisms of schemes.
Write $X = X_{n_0}$, $\sigma = \psi_{n_0} (\phi_{n_0})^{-1}$,
and let $\mc{L} = i_{n_0}^*\mc{O}(1,0,0,\dots, 0)$.
Then
\begin{enumerate}
\item
There is a ring homomorphism $\varphi: A \to B = B(X, \mc{L}, \sigma)$.
\item
Let $J$ be the ideal of $A$ defined as
\[
J = \{x \in A | Mx = 0\ \text{for all}\ R\text{-point modules}
\ M, \text{all commutative algebras}\ R \}.
\]
Then $J$ is equal to $\ker \varphi$ in all degrees $\geq n_0$.
\item
The sheaf $\mc{L}$ is $\sigma$-ample on $X$, and $\varphi$ is
surjective in large degree.
\item The homomorphism in part (1) is canonically determined by the ring $A$.
\end{enumerate}
\end{theorem}

\begin{proof}
(1) We have the homomorphism of graded rings
$$\theta: A \cong T/I \to B' = \bigoplus_{n =0}^{\infty}
\HB^0(X_n, \mc{M}_n)$$
given by Lemma~\ref{xxlem4.1}(1). With our assumption that
the maps $\phi_n$ and $\psi_n$ are isomorphisms for $n \geq n_0$,
we will show there is also a ring homomorphism
\[
\rho: B' = \bigoplus_{n = 0}^{\infty} \HB^0(X_n, \mc{M}_n)
\to B = \bigoplus_{n = 0}^{\infty} \HB^0(X, \mc{L}_n) =
B(X, \mc{L}, \sigma)
\]
where $\mc{L}_n = \mc{L} \otimes \mc{L}^{\sigma} \otimes
\dots \otimes \mc{L}^{\sigma^{n-1}}$ on $X$.  
Let $n \geq 0$ be arbitrary.  If $n > n_0$ then we set $\tau = (\phi_{n,n_0})^{-1}$, 
which makes sense since $\phi_{n,n_0} = \phi_{n_0} \phi_{n_0 +1} \dots \phi_n$ 
is an isomorphism; otherwise 
$n \leq n_0$ and we set $\tau = \phi_{n_0, n}$.  In any case $\tau$ is 
a morphism from $X_{n_0}$ to $X_n$.

Now we claim that for each $0 \leq m \leq n-1$ we have 
\begin{equation}
\label{E4.5}
\mc{L}^{\sigma^m} = \tau^* i_n^*\mc{O}(0, \cdots,0,
\overset{m+1}{1}, 0, \cdots,0).
\end{equation}
Let us outline the proof the claim for the case $n = n_0$, so that 
$\tau$ is the identity;  
The proof for other $n$ is very similar and 
is left for the reader.  For convenience, 
we write $\psi = \psi_{n_0}$, $\phi = \phi_{n_0}$, 
$i = i_{n_0}$, and $j = i_{n_0 +1}$. 
For $m = 0$, \eqref{E4.5} is just the definition of $\mc{L}$. 

By definition, we have $\sigma = \psi \phi^{-1}$.  Also, $\phi, \psi: X_{n_0 +1} \to 
X_{n_0}$ are defined as the restrictions of the two maps $\pi_1, \pi_2: 
\mb{P}^{\times n_0 + 1} \to \mb{P}^{\times n_0}$ which project onto the first 
(respectively, last) $n_0$ coordinates.  This implies that $\pi_1 j  = i \phi$ and 
$\pi_2 j = i \psi$. Thus altogether we have 
\begin{align*}
\mc{L}^{\sigma} = \sigma^* \mc{L} & =  (\phi^{-1})^* \psi^* i^*\mc{O}(1,0,\cdots,0) \\
& =  (\phi^{-1})^* j^* \pi_2^* \mc{O}(1,0, \cdots, 0) \\
& =  (\phi^{-1})^* j^* \pi_1^* \mc{O}(0, 1,0, \cdots,0) \\
& =  (\phi^{-1})^* \phi^* i^* \mc{O}(0,1,0, \cdots, 0) \\
& =  i^* \mc{O}(0,1,0,\cdots, 0)
\end{align*}
which verifies \eqref{E4.5} for $m = 1$.  Now \eqref{E4.5} may be proven for arbitrary $m$ 
with $0 \leq m \leq n-1$ by induction.  The proof of the induction step is analogous 
to the preceding calculation and is omitted.

Given \eqref{E4.5}, since tensor product of sheaves commutes
with pullback, we see that 
\[
\mc{L}_n = \mc{L} \otimes \mc{L}^{\sigma} \otimes \cdots
\otimes \mc{L}^{\sigma^{n-1}}
 = \tau^* i_n^*(\mc{O}(1,1, \cdots, 1)) = \tau^* \mc{M}_n
\]
and thus there is a map 
\[
\HB^0(X_n, \mc{M}_n) \overset{\tau^*}{\to} \HB^0(X_{n_0},\tau^*\mc{M}_n) \cong \HB^0(X, \mc{L}_n).
\]
This defines the map $\rho$ in degree $n$.  Note that if $n \geq n_0$, then $\tau$ is an 
isomorphism and so $\rho$ is an isomorphism in degrees $\geq n_0$.
To check that $\rho$ is a ring homomorphism is a formality left to the reader.
Finally, we get the desired morphism $\varphi$
by letting $\varphi = \rho \circ \theta$.

(2)  Let $\varphi = \rho \circ \theta$ be as constructed in part (1).
It is easy to see that the $J$ given in the statement of part (2)
is an ideal.  If $a \in (\ker \theta)_n$, then by
Lemma~\ref{xxlem4.1}(2), $(M_0)a = 0$ for all truncated $R$-point
modules $M$ of length $n+1$ and all commutative algebras $R$.
If $N = \bigoplus_{i =0}^{\infty} N_i$ is any $R$-point module,
then given any $m \geq 0$ the graded piece $N_m$ is a locally free $R$-module, so 
there is a locally free $R$-module $S$ with  
$S\otimes_R N_m \cong R$.  
Then $M = \bigoplus_{i = 0}^{n} S \otimes_R N_{i+m}$
is a truncated $R$-point module of length $n+1$, and so $M_0 a = 0$.  This 
implies that $N_m a = 0$.  Since $m$ was arbitrary, we have 
$Na = 0$ and $a \in J$.  Because the map $\rho$ is an isomorphism
in degrees $\geq n_0$, we see that for $n \geq n_0$, if
$a \in (\ker \varphi)_n$ then $a \in (\ker \theta)_n$ and thus $a \in J$.

Conversely, suppose that $a \in J_n$ with $n \geq n_0$.  Given any
truncated $R$-point module $M$ of length $n+1$, since $\phi_n$ is
an isomorphism for $n \geq n_0$ there exists an $R$-point module
$N = \bigoplus_{i = 0}^{\infty} N_i$ with
$M \cong \bigoplus_{i = 0}^n N_i$.  Since $Na = 0$, $(M_0)a = 0$.
Thus $a \in \ker \theta$ by Lemma~\ref{xxlem4.1}(2)
and so $a \in \ker \varphi$.  We conclude that $J$ and
$\ker \varphi$ agree in degrees $\geq n_0$.

(3) Let $\varphi: A \to B = B(X, \mc{L}, \sigma)$ be the morphism
constructed in part (1) above, and set $A' = \varphi(A)$.
Let $\mc{K}$ be the sheaf of total quotient rings of $X$, and
let $R = \HB^0(\mc{K})$ be the function ring of $X$, with induced
automorphism also called $\sigma$.  Then we may write
$B =  \bigoplus_{n = 0}^{\infty} \HB^0(X, \mc{L}_n) t^n \subset Q
= R[t^{\pm 1}; \sigma]$.  We may also write
$A' = \bigoplus_{n = 0}^{\infty} V_n t^n$, for some
$V_n \subset \HB^0(X, \mc{L}_n)$.

As we saw in part (1), $\mc{L}_{n_0} \cong \mc{M}_{n_0}$.  The sheaf $\mc{M}_{n_0}$ 
is very ample on $X = X_{n_0}$ since it is the restriction
to $X$ of the very ample sheaf $\mc{O}(1,1, \dots, 1)$ on
$\mb{P}^{\times n_0}$.  For $n \geq n_0$, by the construction
of $\varphi$ we have that $V_n$ is exactly the image of the map
\[
\HB^0(\mb{P}^{\times n}, \mc{O}(1,1,\dots, 1)) \overset{i_n^*}{\to}
\HB^0(X_n, \mc{M}_n) \overset{(\phi_{n, n_0})_*}{\to}
\HB^0(X, \mc{L}_n).
\]
Since $\phi_{n, n_0}$ is an isomorphism for $n \geq n_0$, 
it follows that for $n \geq n_0$ the sections in $V_n$ generate the sheaf $\mc{L}_n$ on $X$, and
moreover that the corresponding map $X \to \mb{P}(V_n)$ is a
closed embedding. By Lemma~\ref{xxlem3.4}, $A'$ is a big subalgebra
of $Q$.

Now suppose that $\sigma$ is not quasi-unipotent.  Then by
Proposition~\ref{xxprop3.5} we have that every big finitely
graded subalgebra of $Q$ has exponential growth.  In particular,
then $A'$ has exponential growth, so necessarily $A$ has
exponential growth, contradicting the hypothesis.

Thus $\sigma$ must be quasi-unipotent, and since $\mc{L}_{n_0}$
is very ample on $X$, $\mc{L}$ must be $\sigma$-ample
\cite[Theorem 1.3]{Ke1}.  Now we have all of the needed hypotheses
to apply Lemma~\ref{xxlem4.3} to the injection $A' \to B$.  We
conclude that $(A')_n = B_n$ for $n \gg 0$, in other words
$\varphi$ is a surjection in large degree as required.

(4) Given the ring $A$, we constructed the ``point scheme data'' $X_n, i_n, \phi_n, \psi_n, \mc{M}_n$
using a fixed presentation $A \cong T/I$, where $T = k \langle x_0, \dots, x_d \rangle$ 
and $\dim A_1 = d+1$.  Given a different presentation $A \cong T/I'$, there is some 
linear transformation $T_1 \to T_1$ which induces an automorphism $\delta: T \to T$ such that 
$I' = \delta(I).$  From this it quickly follows that up to isomorphism the point scheme data does 
not depend on the choice of presentation.

Thus the map $\varphi$ we constructed above depends only on the given ring $A$ and possibly 
the choice of integer $n_0$.  If we choose some different $n_1$ instead, such that $\phi_n$ and $\psi_n$ are isomorphisms 
for all $n \geq n_1$, then the construction produces 
alternative data $X', \mc{L}', \sigma'$ and a map $\varphi': A \to B(X', \mc{L}', \sigma')$.  
If $n_1 > n_0$ then set $\tau = \phi_{n_1, n_0}$, while if $n_1 < n_0$ then set $\tau = (\phi_{n_0, n_1})^{-1}$.
In either case $\tau: X' \to X$ is an isomorphism, $\mc{L}' = \tau^*(\mc{L})$, and 
$\sigma' = \tau^{-1} \sigma \tau$.  Thus pullback of sections via $\tau$ defines an isomorphism 
$\gamma: B(X, \mc{L}, \sigma) \to B(X', \mc{L}', \sigma')$ such that 
$\gamma \varphi = \varphi'$.  Thus the map $\varphi$ is canonically determined.
\end{proof}

We can now derive Theorem~\ref{xxthm1.1} as a special case
of the result just proved.

\begin{proof}[Proof of Theorem~\ref{xxthm1.1}]
Assume the notation from the beginning of this section,
and let $A$ be a strongly noetherian graded algebra
generated in degree $1$. The Hilbert
scheme theorem \cite[Corollary E4.5]{AZ2} shows that since
$A$ is strongly noetherian, the morphisms of schemes
$\phi_n: X_{n+1} \to X_n$ are isomorphisms for all $n \gg 0$.
As is shown in \cite[Proposition 10.2]{KRS}, by considering
left instead of right modules, the strong noetherian property
for $A$ also implies that the maps $\psi_n: X_{n+1} \to X_n$
are isomorphisms for all $n \gg 0$. Then we may choose $n_0$
such that $\phi_n$ and $\psi_n$ are isomorphisms for all
$n \geq n_0$, and then $X = X_{n_0}$ is the point scheme for $A$.
Since $A$ is noetherian, \cite[Theorem 0.1]{SZ} shows that
$A$ has subexponential growth.  Now all of the hypotheses of
Theorem~\ref{xxthm4.4} are satisfied, and so there exists a ring
homomorphism $\varphi: A \to B(X, \mc{L}, \sigma)$.  By Theorem~\ref{xxthm4.4}(4), 
this homomorphism is canoncially determined.
\end{proof}

We now prove several consequences of the main theorem. See
\cite{AZ1} for the definition of a noncommutative projective
scheme and the definitions of ``cohomological dimension''
and ``$\chi$-condition''.
\begin{corollary}
\label{xxcor4.6}
Let $A$ be a strongly noetherian connected graded algebra generated in degree 1.
Let $\varphi: A \to B(X, \mc{L}, \sigma)$ be the map given by
Theorem~\ref{xxthm1.1}.
\begin{enumerate}
\item Suppose that $A \subset Q = R[t^{\pm 1}; \tau]$, where $R$ is a
commutative algebra, and that there is a unit $u \in R$ with $ut \in A_1$. 
Then the map $\varphi$ is an isomorphism in large degree. In particular,
there is an isomorphism of noncommutative projective schemes
$\proj A \cong (\coh X, \mc{O}_X)$.
\item
Suppose that $A$ has a faithful $k$-point module $M$. If 
$A^{(n)}$ is prime for all $n\geq 2$,
Then $A$ is a domain.
\item
Let $M$ be a point module for $A$, and let $A'= A/\rann(M)$.
Then $A'$ has finite GK-dimension, finite Krull dimension,
finite cohomological dimension, and satisfies the
$\chi$-condition.
\end{enumerate}
\end{corollary}

\begin{proof}
(1)
Since $R$ is commutative, we may define 
an $A\otimes R$-module structure on $Q$ by the rule 
$$q\cdot (a\otimes r)=rqa$$
for all $q\in Q$, $a\in A$ and $r\in R$.  

For all $n \geq 0$ 
we have $Q_n (A_1) \supset R t^n u t =
R u^{\sigma^n} t^{n+1} = R t^{n+1} = Q_{n+1}$.  Thus $M := Q_{\geq 0}$ is a cyclic 
$A\otimes R$-module and so is an $R$-point module for $A$.
Since $A$ embeds in $M$, we have $\rann_A M = 0$.  By Theorem~\ref{xxthm1.1}(1),
$(\ker \varphi)_n = 0$ for $n \gg 0$; by Theorem~\ref{xxthm1.1}(2) $\varphi$ is 
also surjective in large degree.  Then $A$ and $B$ are
isomorphic in large degree and $\proj A \cong \proj B$.
Further, it is known that $\proj B \cong (\coh X, \mc{O}_X)$
\cite[Theorem 3.12]{AV}.

(2) If $M$ is a faithful point module for $A$, then it is critical of GK-dimension 1. 
Hence $A$ is graded prime, hence prime as an ungraded ring. 
Since $A$ has a faithful point module, the map $\varphi: A \to B:=B(X,\mc{L},\sigma)$ 
is bijective in large degree by Theorem \ref{xxthm1.1}.  Since $A$ is prime, 
$\varphi$ is injective in all degrees. It remains to show that $B$ is a domain.  For 
this it suffices to show that $X$ is integral.

Since the Veronese ring 
$A^{(n)}$ is prime, so is $B^{(n)}$ for all $n\geq 2$. Note that $B^{(n)} \cong 
B(X,\mc{L}_n,\sigma^n)$.  Without changing $X$, we may replace $B$ by $B^{(n)}$ and 
$\sigma$ by $\sigma^n$.  By doing so we may assume that $\sigma$ fixes all of the 
irreducible components of $X$.  Since $B$ is prime, $X$ has only one irreducible 
component by \cite[Theorem 4.4]{AS}.  It follows by another application of 
\cite[Theorem 4.4]{AS} that $X$ is reduced, so it is integral.

(3) Replacing $A$ by $A'$ we may assume that $A$ has a faithful point module. Then by 
Theorem \ref{xxthm1.1} the map $\varphi: A\to B:=B(X,\mc{L},\sigma)$ is bijective in 
large degree.  It is known that $B$ has finite GK-dimension, finite Krull dimension, 
finite cohomological dimension, and satisfies the $\chi$-condition 
(see \cite[Theorem 1.2]{Ke1}, \cite[Thoerem 3.12]{AV}, and 
\cite[Theorem 4.5, Theorem 7.4]{AZ1}).  Then $A$ also 
has all of these properties.
\end{proof}

\begin{proof}[Proof of Theorem~\ref{main cor}.]
Let $A$ be semiprime, strongly noetherian, connected graded and generated in degree $1$, 
with graded ring of fractions of the form 
$Q = R[t^{\pm 1}; \tau]$ where $R$ is commutative.  Then by the graded Goldie's theorem 
$A$ has a homogeneous nonzerodivisor of some degree $n \geq 1$, so $Q A_n = Q$.
Since $Q$ is graded artinian and $A_n=(A_1)^n$, it follows that $Q A_1 =Q$.  
This shows that $Q_{\geq 0}$ is generated by $Q_0$ as an $A \otimes R$-module 
and hence is an $R$-point module for $A$.  Now just as in the proof of Corollary~\ref{xxcor4.6}(1), 
this implies that the canonical map $\varphi$ is an isomorphism in large degree.
\end{proof}

\begin{remark}
\label{xxrem4.7} It would be interesting to know what happens when the hypotheses of 
Theorem~\ref{xxthm1.1} and Theorem~\ref{main cor} are weakened.  Suppose we start 
with a ring $A$ which is generated in degree $1$ but is only noetherian rather than 
strongly noetherian.  In this case a projective point scheme $X$ may not exist 
\cite{Ro, KRS}, and there is no hope of producing a canonical homomorphism to a 
twisted ring by the same method.  It is possible there may be an analog of some parts 
(e.g., part (2)) of Corollary~\ref{xxcor4.6} for noetherian but not strongly 
noetherian rings, but any classification result along the lines of Theorem~\ref{main cor} 
will necessarily include more exotic examples such as the na{\"i}ve blowups defined in 
\cite{KRS}.  One might also try to relax the hypothesis of generation in degree $1$ in 
Theorem~\ref{xxthm1.1}, while keeping the strong noetherian hypothesis.  In this case 
the point scheme $X$ exists but it is not clear how to define the automorphism 
$\sigma$ in general.  Still, in the special case of graded rings of GK-dimension 2 
which are strongly noetherian but not generated in degree $1$, it is known that 
(possibly after replacing by a Veronese subring) such rings are idealizer rings 
inside twisted rings \cite{AS}.  So there is hope that some homomorphism to a twisted 
ring may still be defined in general regardless of generation in degree $1$.
\end{remark}

\section{Rings $Q$ of finite GK-type}
\label{sect5}

Let $R$ be a commutative algebra with automorphism $\sigma$ and
let $Q$ be the skew Laurent polynomial ring $R[t^{\pm 1}; \sigma]$.
Suppose $Q$ is of finite GK-type (see Definition \ref{xxdef2.1}).
Question \ref{xxque1.4} asks:
\begin{enumerate}
\item
Does every big finitely graded subring of $Q$ have the
\emph{same} GK-dimension?
\item
Must the algebra $Q$ itself has finite GK-dimension?
\end{enumerate}
Our aim in this last section is to give positive answers to both
of these questions in a special case, which was stated in
Theorem~\ref{xxthm1.5}.  Throughout this section we will restrict
our attention to the following situation.  
\begin{hypothesis}
\label{xxhyp5.1}
Let $k$ be an algebraically closed field, and 
let $X$ be an integral projective $k$-scheme with
quasi-unipotent automorphism $\sigma$.  Let $\mc{K}$ be the
sheaf of total quotient rings of $X$, which is the constant
sheaf on the function field $K = k(X)$.  Let $Q =
K[t^{\pm 1}; \sigma]$, where the automorphism $\sigma$ of $K$
is induced by the automorphism $\sigma$ of $X$.
\end{hypothesis}

We claim that under Hypothesis~\ref{xxhyp5.1} the ring $Q$ is
of finite GK-type. Indeed, if $\mc{L}$ is a very ample invertible
sheaf on $X$, then $\mc{L}$ is $\sigma$-ample
\cite[Theorem 1.3]{Ke1} and the ring $B = B(X, \mc{L}, \sigma)$
has finite GK-dimension \cite[Theorem 6.1]{Ke1}.  Since the ring
$B$ is big in $Q$ by Lemma~\ref{xxlem3.4}, $Q$ is of finite
GK-type by Proposition~\ref{xxprop1.3}. 

We note that given an automorphism $\sigma$ of
a field $K$ of finite transcendence degree over $k$, it is possible 
that there does not exist any projective model $X$ of $K$ with a regular 
automorphism $\tau \in \aut(X)$ inducing the map $\sigma$ on rational functions.  
In this case, we call $\sigma$ \emph{non-geometric}.
We have verified Question \ref{xxque1.4}(2) for some specific 
examples in the non-geometric case, but the details are omitted here. In general
Question \ref{xxque1.4} is open for non-geometric $\sigma$.

We start with some easy lemmas and assume Hypothesis~\ref{xxhyp5.1}
throughout.

\begin{lemma}
\label{xxlem5.2}
Let $V$ be a finite-dimensional $k$-subspace of $K$.  Then there is some
very ample invertible sheaf $\mc{L} \subset \mc{K}$ such that
$V \subset \HB^0(\mc{L}) \subset K$.
\end{lemma}

\begin{proof}
Let $\mc{M} \subset \mc{K}$ be any very ample invertible sheaf
on $X$, and let $W = \HB^0(\mc{M})$.  By Lemma~\ref{xxlem3.4},
we can choose $u \in W$ such that $K$ is the field of fractions
of $k[U]$, where $1 \in U = u^{-1}W$.  Putting the elements of
$V$ over a common denominator $z$, we have $V \subset z^{-1}
U^n$ for some $n \geq 0$.  Then $\mc{L} = z^{-1}u^{-n}
\HB^0(\mc{M}^{\otimes n}) \mc{O}_X$ is a very ample sheaf
(isomorphic to $\mc{M}^{\otimes n}$) such that $V
\subset z^{-1}u^{-n} W^n \subset \HB^0(\mc{L})$.
\end{proof}

\begin{lemma}
\label{xxlem5.3}
Let $A = \bigoplus_{n = 0}^{\infty} V_n t^n$ be a big
finitely graded subalgebra of $Q$. Let $W$ be a
finite-dimensional $k$-subspace of $K$. Then there is $0
\neq u \in K$ and $n \geq 1$ such that $Wu \subseteq V_n$.
\end{lemma}

\begin{proof} Since $Q$ is of finite GK-type, $A$ has finite
GK-dimension. Since $A$ is a domain, the set of all nonzero homogeneous
elements of $A$ is an Ore set \cite[Proposition 4.13]{KL}. Since
$A$ is big in $Q$, the graded ring of fractions of $A$ must be
$K[(t^m)^{\pm 1}; \sigma^m]$ for some $m$ by Lemma~\ref{xxlem2.3}(2).
Now using the Ore condition for $A$, we may write the elements in
$W$ over a common denominator.  So for some $n \geq 1$ and $z
\in A_n$, $W \subseteq A_n z^{-1}$.  This means that $W u
\subseteq V_n$ for $u = z t^{-n} \in K$.
\end{proof}

\begin{lemma}
\label{xxlem5.4}
Let $\mc{L}$ be a very ample invertible sheaf on $X$.
If $A$ is the subring of $B = B(X, \mc{L}, \sigma)$ generated
by $B_1$, then $A_n = B_n$ for all $n \gg 0$.
\end{lemma}

\begin{proof}
Fix some embedding $\mc{L} \subset \mc{K}$ and write $B =
\bigoplus_{n = 0}^{\infty} \HB^0(X, \mc{L}_n) t^n \subset Q$.
Let $A = \bigoplus_{n = 0}^{\infty} V_n t^n$.  Since $\mc{L}$ is
very ample, it is generated by its global sections $V_1$, and
determines an embedding $X \to \mb{P}(V_1)$.  Then for each
$n \geq 1$, $\mc{L}_n$ is generated by $V_n = V_1 V_1^{\sigma}
\cdots V_1^{\sigma^{n-1}}$, and the sections $V_n$ determine
an embedding $X \to \mb{P}(V_n)$.  Since $\sigma$ is quasi-unipotent,
$\mc{L}$ is $\sigma$-ample \cite[Theorem 1.3]{Ke1}, so
all of the hypotheses of Lemma~\ref{xxlem4.3} are
verified and that result gives $A_n = B_n$ for $n \gg 0$.
\end{proof}

Now we see how to compare an arbitrary big subalgebra of $Q$ with twisted
rings inside $Q$.

\begin{proposition}
\label{xxprop5.5}
Let $A$ be a big connected finitely graded subring of $Q$.
\begin{enumerate}
\item There is a very ample sheaf $\mc{L} \subset \mc{K}$
such that $A \subset B(X, \mc{L}, \sigma) \subset Q$.
\item There is a very ample sheaf $\mc{M}$ and $h \geq 1$ such that
$B(X, \mc{M}, \sigma^h) \subset A^{(h)}$.
\end{enumerate}
\end{proposition}

\begin{proof}
(1) Write $A = \bigoplus_{n = 0}^{\infty} V_n t^n$.
If $A$ is generated by elements of degree less than $n$, Then $A$
is contained in the connected graded ring generated in degree $1$ by
$Vt$, where $V=k+\sum_{i=1}^n V_i$.  Now by Lemma~\ref{xxlem5.2}, we
can find a very ample invertible sheaf $\mc{L} \subset \mc{K}$ such
that $V \subseteq \HB^0(\mc{L}) = W$.  Then
\[
A \subseteq k[Vt]
\subseteq k[Wt] \subseteq B(X, \mc{L}, \sigma) =
\bigoplus_{n = 0}^{\infty} \HB^0(\mc{L}_n) t^n \subseteq Q.
\]

(2) Pick an arbitrary very ample invertible sheaf $\mc{L}' \subset
\mc{K}$, and let $W' = \HB^0(\mc{L}') \subset K$.  By
Lemma~\ref{xxlem5.3}, we have $W' u \subseteq V_n$ for some
$n \geq 1$ and some $0 \neq u \in K$.  Write $W = W'u$.
Then let $\mc{L} = W \mc{O}_X \subset \mc{K}$; so $\mc{L}$ is a
very ample sheaf isomorphic to $\mc{L}'$ but embedded differently
in $\mc{K}$.

Let $C$ be the subring of $B = B(X, \mc{L}, \sigma^n) \subset
K[s; \sigma^n]$ generated by $B_1$.  Then
\[
C_m = W W^{\sigma^n} \dots W^{\sigma^{(m-1)n}}s^m \subset
V_n V_n^{\sigma^n} \dots V_n^{\sigma^{(m-1)n}} s^m \subset V_{nm} s^m
\]
for all $m$ and so we have an inclusion of connected graded rings
$C \subset A^{(n)} = \bigoplus_{m \geq 0} A_{mn}$.  Since $\sigma$
is quasi-unipotent by hypothesis, the automorphism $\sigma^n$ is
also quasi-unipotent.  Then by Lemma~\ref{xxlem5.4}, $C$ and $B$
are equal in large degree, say in degrees $\geq p$.  Therefore
passing to a further Veronese subring we have an inclusion
$B^{(p)} = B(X, \mc{M}, \sigma^{pn}) \subset A^{(np)}$,
where $\mc{M} = \mc{L} \otimes \mc{L}^{\sigma^n} \otimes
\dots \otimes \mc{L}^{\sigma^{n(p-1)}}$ is very ample.
\end{proof}

The following is a result of Keeler \cite{Ke1}.
\begin{lemma}
\label{xxlem5.6}
The number $\GK B(X,\mc{N},\sigma^n)$ is the same positive integer for
every ample invertible sheaf $\mc{N}$ and all $n \geq 1$.
\end{lemma}
\begin{proof} Since by Hypothesis~\ref{xxhyp5.1} $\sigma$ is
quasi-unipotent, this follows from \cite[\S~6]{Ke1}: Theorem 6.1,
Equation 6.3, Lemma 6.5(2) and Lemma 6.8.
\end{proof}

We are now ready to prove the main result of this section.

\begin{proof}[Proof of Theorem~\ref{xxthm1.5}]
Note that the hypotheses of the theorem allow us to assume
Hypothesis~\ref{xxhyp5.1}. Let $d=\GK B(X,\mc{L},\sigma)$.

(1) Let $A$ be any big finitely graded subalgebra of $Q$.  
Since the ring
$k + A_{\geq 1}$ has the same GK-dimension, we may assume that $A$ is connected.
Since $A^{(h)} \subset A$ as ungraded rings,
by Proposition~\ref{xxprop5.5}
for some $h$ and very ample sheaves $\mc{N}, \mc{M}$ we have
\[
\GK B(X, \mc{M}, \sigma^h) \leq \GK A^{(h)} \leq \GK A \leq
\GK B(X, \mc{N}, \sigma),
\]
where the two ends are equal to $d$ by Lemma \ref{xxlem5.6}.
Thus $\GK A = d$.

(2) Clearly $\GK Q\geq \GK A=d$. It remains to show that
$\GK Q\leq d+\dim X$.

Let $V$ be a finite-dimensional $k$-subspace of $Q$.  There is
a finite-dimensional $W \subseteq K$ and $m \geq 0$ such that
$V \subseteq V': = \sum_{i = -m}^m W t^i$; we may choose $W$ so
that $1 \in W$, and so that $W = \HB^0(\mc{L})$ for
some very ample invertible sheaf $\mc{L} \subset \mc{K}$, by
Lemma~\ref{xxlem5.2}.  We want to give an upper bound for $\GK k [V']$.
Setting $V'' = Wt^{-1} + W + Wt$, then $(V'')^n \supset V'$,
so $\GK k [V'] \leq \GK k [V'']$ and it is enough to bound the
latter dimension.

By the fact that $1\in W$ we have
\[
(V'')^n \subseteq \sum_{i = -n}^{n}
{(W^{\sigma^{-n}} W^{\sigma^{-n+1}} \dots W^{\sigma^{n}})}^n t^i.
\]
Note also that
\[
{(W^{\sigma^{-n}} W^{\sigma^{-n+1}} \cdots W^{\sigma^{n}})}^n
\subseteq \HB^0({(\mc{L}^{\sigma^{-n}} \otimes \mc{L}^{\sigma^{-n+1}}
\otimes \cdots \otimes \mc{L}^{\sigma^{n}})}^{\otimes n}) \cong
\HB^0(\mc{L}_{2n+1}^{\otimes n})
\]
So altogether we may conclude that
\begin{equation}
\label{E5.7}
\dim_k (V'')^n \leq (2n+1) \dim_k \HB^0(\mc{L}_{2n +1}^{\otimes n}).
\end{equation}

Now for each $n \geq 0$ write $e_n = \dim_k \HB^0(\mc{L}_n)$,
$f_n = (2n+1) \dim_k \HB^0(\mc{L}_{2n +1}^{\otimes n})$,
and $g_n = \dim_k \HB^0(\mc{L}_n^{\otimes n})$. Since
$B =B(X, \mc{L},\sigma)$ is a finitely generated domain, we have
$d = \GK B = (\overline{\lim} \log_n e_n) + 1$.  It follows from
\eqref{E5.7} that  $\GK k [V''] \leq \overline{\lim} \log_n f_n$.
Since $\mc{L}$ is very ample, we have $\dim_k
\HB^0(\mc{L}_{2n +1}^{\otimes n}) <
\dim_k \HB^0(\mc{L}_{2n+1}^{\otimes 2n+1})$, and so
$\overline{\lim} \log_n f_n = (\overline{\lim} \log_n g_n) +1$.
So we will be done if we can show that
$g_n \leq n^{\dim X} e_n$ for $n \gg 0$.

Since $\mc{L}$ is very ample, for $n \gg 0$ we have
$\HB^q(\mc{L}^{\otimes n} \otimes \mc{M})= 0$ for any ample
sheaf $\mc{M}$ and all $q > 0$ \cite[p. 520, Theorem 1]{Fj}.
This implies that, for $n \gg 0$, $\HB^q(\mc{L}_n^{\otimes n}) = 0$.
Also, $\HB^q(\mc{L}_n) = 0$ for $n \gg 0$ and all $q > 0$ since
$\mc{L}$ is $\sigma$-ample by \cite[Theorem 1.3]{Ke1}.
For a coherent sheaf $\mc{F}$ on $X$, define
$\chi(\mc{F}) = \bigoplus_{q = 0}^{\infty} (-1)^q \dim_k \HB^q(\mc{F})$.
So $g_n = \chi(\mc{L}_n^{\otimes n})$ and $e_n= \chi(\mc{L}_n)$ for
$n \gg 0$.

Now, as in \cite[p. 529]{Ke1}, the Riemann-Roch formula
may be used to show for an invertible sheaf $\mc{M}$ that
\[
\chi(\mc{M}) = \sum_{j = 0}^{\dim X} \sum_{V \in A_j}
a_{V} \overset{j}{\overbrace{(\mc{M} . \mc{M} . \cdots . \mc{M})}}_V,
\]
where $A_j$ is some finite set of $j$-dimensional subvarieties of
$X$ and the $a_V$ are some rational numbers. All that we need to
know about the intersection form appearing in this formula is that
it is multilinear (with respect to addition of divisors, which
corresponds to tensor product of sheaves).  So we conclude for
all $n \geq 0$ that
\[
\chi(\mc{L}_n^{\otimes n}) = \sum_{j = 0}^{\dim X} n^j
\sum_{V \in A_j} a_{V} ((\mc{L}_n)^j)_V
\leq n^{\dim X} \chi(\mc{L}_n)
\]
and so $g_n \leq n^{\dim X} e_n$ for $n \gg 0$ as we wished.
\end{proof}

The lower bound in Theorem~\ref{xxthm1.5}(2) holds when
$\sigma = Id$ (or $\sigma$ has finite order).  
We can show that the 
upper bound is achieved in some specific examples, for instance if
$X = \mb{P}^n$ and $\sigma$ is a generic automorphism.
In general it is not clear exactly how the value of $\GK Q$
depends on the properties of the automorphism $\sigma$.

To close, we show that as a corollary of Theorem~\ref{xxthm1.5}
we can obtain the exact value of the Gelfand-Kirillov
transcendence degree of the division rings of fractions of
twisted rings. The \emph{GK-transcendence degree} of an
algebra $A$ is defined to be
$$\GKtr A=\sup_{V} (\inf_{z} \GK k[zV] )$$
where $V$ ranges over all finite-dimensional subspaces
of $A$ containing $1$ and $z$ ranges over all
nonzerodivisors of $A$.  This is a very useful invariant for
division rings, but it is difficult to calculate in general
\cite{GK,Zh}.

\begin{corollary}
\label{xxcor5.8}
Let $B = B(X, \mc{L}, \sigma)$ where $\mc{L}$ is $\sigma$-ample,
and let $D$ be the division ring of fractions of the noetherian
domain $B$. If $A$ is any algebra with $B\subset A \subset D$,
then $\GKtr A=\GK B$.
\end{corollary}

\begin{proof} (1) By \cite[Propositions 2.1 and 3.1(3)]{Zh},
\[
\GKtr D \leq \GKtr A \leq \GKtr B \leq \GK B.
\]
It remains to show that $\GKtr D \geq \GK B$. Let $\nu: B \to B$
be the valuation defined by the grading of $B$.  By
\cite[Proposition 6.5]{Zh},  $\nu$ can be extended to a valuation
$\nu: D\to Q$ where $Q = K[t^{\pm 1}; \sigma]$ is the graded ring
of fractions of $B$. Let $V$ be a finite-dimensional subspace of
$K$ containing $1$ such that $K$ is the field of fractions of
$k[V]$.  Set $W=Vt+V+Vt^{-1}$. It will suffice to show that
$\GK k[zW] \geq \GK B$ for all nonzero $z \in Q$. By
\cite[Theorem 6.7]{Zh}, we only need to show that
$\GK k[\nu(zW)] = \GK k[\nu(z) W] \geq \GK B$.  Thus
we can assume that $z$ is a homogeneous element of $Q$, and we write
$z=q t^n$ for some $q \in K$ and $n \in \mb{Z}$.  If $n \geq 0$, then 
$k[zW]$ contains the connected $\mb{N}$-graded algebra $C= k[qVt^{n+1}]$, 
which is big in $Q$ by definition.  Then $\GK C = \GK B$ by Theorem~\ref{xxthm1.5}.
Otherwise $n < 0$, and $k[zW]$ contains the algebra $D = k[qVt^n]$.  The ring $D$ has 
the same GK-dimension as the big $\mb{N}$-graded algebra $C = k[qVt^{-n}]$, and again 
$\GK C = \GK B$.  In either case we have that $\GK k[zW] \geq \GK B$.
\end{proof}

\section*{acknowledgments}
We thank Mike Artin, Dennis Keeler, 
Zinovy Reichstein, Paul Smith, and Toby Stafford
for helpful conversations.

\providecommand{\bysame}{\leavevmode\hbox to3em{\hrulefill}\thinspace}
\providecommand{\MR}{\relax\ifhmode\unskip\space\fi MR }
\providecommand{\MRhref}[2]{
  \href{http://www.ams.org/mathscinet-getitem?mr=#1}{#2}
}
\providecommand{\href}[2]{#2}

\end{document}